# Formal Languages and Infinite Groups

ROBERT H. GILMAN

ABSTRACT. This article is an introduction to formal languages from the point of view of combinatorial group theory.

## Contents

## 1. Introduction

A formal language is a subset of a finitely generated free monoid. Certain classes of formal languages play important roles in computer science and are

1991 *Mathematics Subject Classification.* Primary 20F32; Secondary 20F05 68Q42 68Q45.
The author was supported in part by NSF Grant DMS-9401090.
This paper is in final form and no version of it will be submitted for publication elsewhere.

beginning to have interesting applications to combinatorial group theory. Automatic groups for example are defined in terms of rational languages and rational relations. Rational languages also occur in the study of word hyperbolic groups and are useful in computational group theory. The larger class of context-free languages can be used to characterize virtually free and virtually cyclic groups. The fundamental groups of compact geometrizable three-manifolds are describable in terms of indexed languages.

This expository article is an introduction to formal language theory for people whose interests include combinatorial group theory. We concentrate on the basic properties of rational and context-free languages because these are the languages which have had most application to group theory so far. Rational languages are easily described; they are the closure of the finite subsets of finitely generated free monoids under union, product and generation of submonoids. Rational subsets of arbitrary monoids are defined in the same way, but they do not possess all the properties of rational languages.

The rational subsets of a monoid $M$ are precisely the subsets accepted by finite automata over $M$. A finite automaton $\Gamma$ over $M$ is a finite directed graph with a distinguished initial vertex, some distinguished terminal vertices, and with edges labelled by elements from $M$. The set accepted by $\Gamma$ is the collection of labels of paths from the initial vertex to a terminal vertex. Figure 1 shows a finite automaton over $\Sigma^*$, the free monoid over a finite set $\Sigma$ containing $a$ and $b$. The symbol $\epsilon$ stands for the unit element of $\Sigma^*$. In this figure and in all others an arrow with no source indicates the initial vertex, and arrows with no target identify terminal vertices. This finite automaton accepts the rational language $\{a^i b^j \mid i \geq 0, j \geq 0\}$.

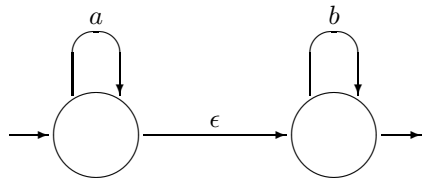

FIGURE 1. A finite automaton over $\Sigma^*$.

There are also acceptors, called pushdown automata, for context-free languages. Since the standard definition of pushdown automaton although useful for understanding aspects of programming is unwieldy, we do not use it. For us a pushdown automaton accepting a context-free language in $\Sigma^*$ is just a finite automaton over $M_{cf} \times \Sigma^*$ where $M_{cf}$ is a certain monoid characterizing context-free languages. A word $w \in \Sigma^*$ is accepted by such a pushdown automaton if there is a path from initial vertex to terminal vertex with label $(1, w)$. This approach follows that of several authors. See Section 10 for references.

Figure 2 shows a pushdown automaton. $P_d$ and $Q_d$ are elements of $M_{cf}$, and as before $a, b \in \Sigma$. This finite automaton accepts the rational subset $\{(P_d^i Q_d^j, a^i b^j) \mid$ $i \geq 0, j \geq 0\}$ of $M_{cf} \times \Sigma^*$. Construed as a pushdown automaton it accepts the set of all words $a^i b^j$ such that $P_d^i Q_d^j = 1$ in $M_{cf}$. Given that $P_d Q_d = 1$ and $P_d$ and $Q_d$ have infinite order, we see that this pushdown automaton accepts the language $\{a^n b^n \mid n \geq 0\}$.

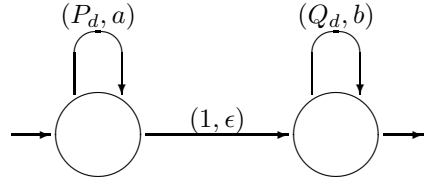

FIGURE 2. A finite automaton over $M_{cf} \times \Sigma^*$.

Replacing $M_{cf}$ by different monoids affords a uniform way of defining automata accepting other classes of languages. More precisely we specify $M$ and a subset $X \subset M$, and consider automata over $M \times \Sigma^*$ which accept $w \in \Sigma^*$ if for some $x \in X$ there is a path with label $(x, w)$ from initial vertex to terminal vertex. Theorem 6.2 characterizes the classes which can be defined this way; they fall between full trios and full abstract families of languages. Stack automata are considered briefly from this point of view in Section 8.

We make the following conventions. $\Sigma^*$ denotes the free monoid over a set $\Sigma$. $\Sigma$ may be referred to as an alphabet and elements $w \in \Sigma^*$, as words. The unit element of $\Sigma^*$ is $\epsilon$, the empty word; and $|w|$ is the length of the word $w$. A language contained in $\Sigma^*$ is also called a language over $\Sigma$. $\Sigma$ and $\Delta$ will always denote nonempty finite sets. In order to avoid some awkward constructions we write function and relation symbols on the right.

## 2. Rational Sets

DEFINITION 2.1. The rational subsets of a monoid $M$ are the closure of its finite subsets under union, product, and monoid closure.

We call union, product and monoid closure the rational operations and write them as $S_1 \cup S_2$, $S_1 S_2$, and $S_1^*$ respectively. Notice that this notation is consistent that for free monoid adopted above. Rational subsets of free monoids are called rational languages and also regular languages.

THEOREM 2.2. *Let $M$ and $M'$ be monoids and $\sigma : M \to M'$ a homomorphism.*
  (i) *For any set of generators $A$ of $M$, each rational subset of $M$ lies the submonoid generated by a finite subset of $A$.*
 (ii) *The rational subsets of $M$ are closed under rational operations.*
(iii) *The image of a rational subset of $M$ under $\sigma$ is a rational subset of $M'$; and if $\sigma$ is onto, every rational subset of $M'$ is obtained this way.*

PROOF. We prove the first assertion by induction on rational sets. Let $S$ be a rational subset of $M$. If $S$ is finite, then clearly it lies the submonoid generated by a finite subset of $A$. Suppose $S = S_1 \cup S_2$ and each $S_i$ lies a submonoid generated by a finite subset of $A$; then $S$ is contained in the submonoid generated by the union of the two finite sets. The same argument works if $S = S_1 S_2$ or $S = S_1^*$, and we conclude that (i) holds.

Part (ii) follows immediately from Definition 2.1, and the last part is proved by induction on rational sets. ☐

Finite automata afford a convenient way of describing individual rational sets.

DEFINITION 2.3. A finite automaton defined over a monoid $M$ is a finite directed graph with
  (i) Edges labelled by elements of $M$;
  (ii) An initial vertex;
  (iii) A set of terminal vertices.
A finite automaton over a subset $A \subset M$ is just a finite automaton over $M$ whose edge labels are all in $A$.

Since finite automata play such a large role in this exposition, we will refer to them as automata; that is, from now on automaton means finite automaton. By a path in an automaton we mean a directed path. The label of a path is the product in order of the labels of its edges; in particular the label of an empty path, i.e., a path of length zero from a vertex to itself, is the unit element.

DEFINITION 2.4. The set accepted by an automaton is the collection of labels of successful paths. A successful path is one from the initial vertex to a terminal vertex.

Figure 1 shows a finite automaton which accepts the rational language $a^*b^*$.

LEMMA 2.5. *Let $A$ be a subset of a monoid $M$. Every set accepted by an automaton over $A$ is accepted by one such that*
  (i) *The initial vertex has no inedges;*
  (ii) *There a single terminal vertex, and it has no outedges.*

PROOF. Since the set accepted is not changed by the addition of a new isolated terminal vertex, we may assume the set of terminal vertices is not empty. If the initial vertex $p_0$ has inedges, add a new initial vertex $p_0'$ and duplicate only the outedges of $p_0$. For each terminal vertex $p$ with outedges add a new terminal vertex $p'$. Duplicate the inedges of $p$, and remove it from the set of terminal vertices. Changing the first and last edges of a path gives a correspondence between paths from $p_0$ to $p$ and those from $p_0'$ to $p'$. Thus the set of labels of successful paths is not altered by the changes made so far. Since terminal vertices now have no outedges, identifying all terminal vertices does not affect the set of labels of successful paths either. ☐

THEOREM 2.6. *Let $A$ be any set of generators for the monoid $M$. A subset of $M$ is rational if and only if it is accepted by an automaton over $A$.*

PROOF. Argue by induction on rational sets to show that every rational subset of $M$ is accepted by an automaton. Indeed every finite subset of $M$ is accepted by an automaton constructed from a tree, and it is not hard to combine automata satisfying the conditions of Lemma 2.5 to accept unions, products, and monoid closures. If $\Gamma$ accepts $S$ and $\Gamma'$ accepts $S'$, identify the two initial vertices to obtain an automaton accepting $S \cup S'$. Likewise identify the terminal vertex of $\Gamma$ with the initial vertex of $\Gamma'$ to accept $SS'$, and the initial and terminal vertices of $\Gamma$ to accept $S^*$.

For the converse let $\Gamma$ accept $S$ and argue by induction on the number of edges of $\Gamma$. If there are no edges, then $S$ equals $\{\epsilon\}$ or $\phi$ depending on whether or not the initial vertex is terminal. Otherwise construct an automaton $\Gamma_0$ from $\Gamma$ by removing an edge $e$ while retaining its source and target vertices $p$ and $p'$ respectively. Let $a$ be the label of $e$. The label $w$ of any successful path in $\Gamma$ not lying entirely in $\Gamma_0$ is a product

$$(1) \qquad w = w_0 a w_1 a \cdots w_m, \quad m \geq 1$$

for which the following conditions hold. The word $w_0$ is the label of a path in $\Gamma_0$ from $p_0$ to $p$, and each subsequent $w_i$ except the last is the label of a path in $\Gamma_0$ from $p'$ to $p$. Finally $w_m$ is the label of a path in $\Gamma_0$ from $p'$ to a terminal vertex of $\Gamma$. Construct automata $\Gamma_i$, $1 \leq i \leq 3$ from $\Gamma_0$ by altering the initial and terminal vertices of $\Gamma_0$ in the following way.

(i) $\Gamma_1$ has initial vertex $p_0$ and terminal vertex $p$;
(ii) $\Gamma_2$ has initial vertex $p$ and terminal vertex $p'$;
(iii) $\Gamma_3$ has initial vertex $p'$ and the same terminal vertices as $\Gamma$.

For $0 \leq i \leq 3$ let $S_i$ be the set accepted by $\Gamma_i$. These sets are rational by induction, and Equation 1 implies

$$S = S_0 \cup S_1 a (S_2 a)^* S_3.$$

Thus $S$ is rational. $\square$

## 3. Rational Languages

Rational languages satisfy some properties not satisfied by rational sets in general.

By Theorem 2.6 we take automata over free monoids as defined over the free generators unless otherwise stated.

THEOREM 3.1. *For any finite alphabet $\Sigma$*
  (i) *The rational languages over $\Sigma$ are the inverse images of subsets of finite monoids $F$ under homomorphisms $\sigma : \Sigma^* \to F$; and*
  (ii) *The rational languages over $\Sigma$ are closed under Boolean operations.*

PROOF. The second assertion follows immediately from the first. Let $L$ be a rational language contained in $\Sigma^*$ and assume $L$ is accepted by a finite automaton $\Gamma$ over $\Sigma$. For each $w \in \Sigma^*$ define a binary relation $\sim_w$ on the vertices of $\Gamma$ by $p \sim_w p'$ if and only if there is a path from $p$ to $p'$ with label $w$. Because $\Sigma^*$ is free, $\sim_{wv}$ equals the composite $\sim_w \circ \sim_v$. Hence the map $w \to \sim_w$ is a homomorphism from $\Sigma^*$ to the finite monoid of binary relations on the vertices of $\Gamma$. $L$ is the inverse image of the set of all relations $\sim$ such that $p_0 \sim p$ for some terminal vertex $p$.

Conversely suppose that $\sigma : \Sigma^* \to F$ is a homomorphism to a finite monoid $F$, and $L = X\sigma^{-1}$ for some $X \subset F$. Let $\Gamma$ be directed graph with vertices $F$ and labelled edges $x \xrightarrow{a} y$ for all $x, y \in F$ and $a \in \Sigma$ such that $x(a\sigma) = y$. Make $\Gamma$ into an automaton with initial vertex 1 and terminal vertices $X$. Since a path from 1 with label $w$ ends at $w\sigma$, $\Gamma$ accepts $L$. □

The automaton constructed in the last part of the proof above has a special property. Every vertex has exactly one outedge labelled by each generator. Such automata are called *deterministic*; the others are *nondeterministic*. For automata over monoids which are not free there are other notions of deterministic. It is a fact that each rational language is accepted by a unique smallest deterministic automaton, called the minimal automaton of the language.

COROLLARY 3.2. *Every rational language is accepted by a deterministic automaton.*

THEOREM 3.3. *For every rational language $L$ there is an integer $n$ with the property that if $w \in L$ and $|w| > n$, then $w$ can be written as a product $xyz$ such that $xy^*z \subset L$ and $0 < |y| \leq n$.*

PROOF. Let $L$ be accepted by an automaton over $\Sigma$ with $n$ vertices. Every successful path of length greater than $n$ must have a loop of length strictly between 0 and $n + 1$. □

Theorem 3.3 is called the Pumping Lemma for Rational Languages; its usual first application is to show that the language $\{a^n b^n\}$ is not rational.

## 4. Applications of Rational Languages

Let $G$ be a group. A *choice of generators* for $G$ is a surjective monoid homomorphism $\sigma : \Sigma^* \to G$. We will write $\overline{w}$ for the image, $w\sigma$, of $w \in \Sigma^*$. By Theorem 2.2 every rational subset of $G$ lifts to a rational language over $\Sigma$. Just as group elements are represented by words, rational subsets of groups are represented by rational languages. We will say that a choice of generators $\sigma : \Sigma^* \to G$ has *formal inverses* if $\Sigma$ is a union of pairs $\{a, a^{-1}\}$ and $a^{-1}\sigma = (a\sigma)^{-1}$. We emphasize that $\Sigma$ still generates $\Sigma^*$ freely as a monoid; there is no cancellation in $\Sigma^*$.

The *word problem* of $G$ is the language $W(G) = \{w \mid \overline{w} = 1\}$. $W(G)$ depends on the choice of generators, but the dependence is mild. A precise statement is given in Theorem 6.4. The following theorem characterizes the word problem of finite groups with respect to any choice of generators.

THEOREM 4.1. *The word problem of a finitely generated group $G$ is a rational language if and only if $G$ is finite.*

PROOF. Suppose $\Sigma \to G$ is a choice of generators. If $G$ is finite, then $W(G)$ is rational by Theorem 3.1. For the converse suppose $W(G)$ is accepted by an automaton $\Gamma$. We will show that the order of $G$ is at most equal to the number of vertices in $\Gamma$.

Assume without loss of generality that every vertex of $\Gamma$ is on a successful path. Indeed discarding all edges and vertices not on successful paths does not affect the language accepted by $\Gamma$. For any word $w \in \Sigma^*$ there is a word $v$ representing the inverse of $\overline{w}$. Hence $wv$ is accepted by $\Gamma$, and consequently $w$ is the label of a path starting at the initial vertex $p_0$ of $\Gamma$. In particular every element of $G$ is represented by the label of such a path. To complete the proof it suffices to show that if $w$ and $u$ are labels of path from $p_0$ to the same vertex $p$, then $\overline{w} = \overline{u}$. By our assumption there is a path from $p$ to a terminal vertex. Let the label of this path be $v$. We have that $wv$ and $uv$ both represent the identity in $G$ whence $\overline{w} = \overline{u}$. . ☐

THEOREM 4.2. *The subgroup generated by a rational subset of a group is finitely generated.*

PROOF. Let $\langle S \rangle$ be the subgroup generated by $S \subset G$. Suppose $S$ is rational and is accepted by an automaton $\Gamma$ over $G$. If $S = \phi$, then $\langle S \rangle = \{1\}$, so we may assume $S \neq \phi$. By Lemma 2.5 we may suppose that $\Gamma$ has a single terminal vertex $p_t$; and as above we may discard all vertices and edges not on successful paths. Since $\Gamma$ accepts a nonempty set, $p_t$ is not discarded. It follows from our conditions that we can pick a spanning subtree $\Gamma_0$ of $\Gamma$ with root $p_0$, the initial vertex of $\Gamma$, and with all edges directed away from the root. The label of any successful path in $\Gamma$ is a product

$$g = g_0 h_1 g_1 h_2 g_2 \cdots g_{m-1} h_m g_m.$$

in which each $g_i$ is the label of a path in $\Gamma_0$, and each $h_i$ is the label of an edge $e_i$ not in $\Gamma_0$. Let $x_i$ and $y_i$ be the labels of the paths in $\Gamma_0$ from the root to the source and target vertex of $e_i$ respectively, and let $z$ be the label of the path in $\Gamma_0$ from $p_0$ to $p_t$. Since there is at most one path in $\Gamma_0$ between any two vertices, $x_1 = g_0$. Likewise $x_{i+1} = y_i g_i$ for $1 < i < m$, and $z = y_m g_m$. Thus

$$g = x_1 h_1 y_1^{-1} x_2 h_2 y_2^{-1} \cdots x_m h_m y_m^{-1} z$$

It follows that $S \subset H = \langle z, xhy^{-1} \rangle$, the subgroup generated by $z$ and $xhy^{-1}$ where $h$ ranges over all labels of edges $e$ not in $\Gamma_0$, and $x$ and $y$ are the labels of the paths in $\Gamma_0$ to the source and target of $e$ respectively. To complete the proof it suffices to show that the generators for $H$ lie in $\langle S \rangle$. We have $z \in S$ because $z$ is the label of a successful path. Suppose $e$ is an edge not in $\Gamma_0$, and $xhy^{-1}$ is its corresponding generator. Let $v$ be the label of a path in $\Gamma$ from the target vertex of $e$ to $p_t$. It follows from $xhv, yv \in S$ that $xhy^{-1} = (xhv)(yv)^{-1} \in \langle S \rangle$. ☐

Next we use rational languages to prove that free groups have the Howson property.

THEOREM 4.3. *The intersection of two finitely generated subgroups of a free group is finitely generated.*

PROOF. Let $G$ be a finitely generated free group and $\Sigma \to G$ a choice of free generators with formal inverses. By Theorem 2.2 we may assume that for $i = 1, 2$ the finitely generated subgroup $H_i$ is the image of the rational language $L_i$ over $\Sigma$. If each $L_i$ contains all the freely reduced words which represent elements of $H_i$, then $L_1 \cap L_2$ contains all freely reduced words projecting to $H_1 \cap H_2$. Consequently $L_1 \cap L_2$ projects onto $H_1 \cap H_2$, and $H_1 \cap H_2$ is finitely generated by Theorem 4.2. Thus the following lemma completes the proof. □

LEMMA 4.4. *If $\Sigma = \{a, a^{-1} \ldots\}$ and $L \subset \Sigma^*$ is rational, then the corresponding language of freely reduced words is rational too.*

PROOF. Let $\Gamma$ be an automaton defined over $\Sigma$ and accepting $L$, and let $L_0$ be the corresponding set of freely reduced words. We will make $\Gamma$ into an automaton over $\Sigma \cup \{\epsilon\}$. If for vertices $p$ and $q$, $\Gamma$ has a path from $p$ to vertex $q$ whose label freely reduces to the identity, add an edge from $p$ to $q$ with label $\epsilon$ unless such an edge is not already present. This effect of this addition is to add to $L$ some free reductions of words already there. Continue the edge addition procedure as long as possible. Since $\Gamma$ is finite, the process must terminate; and it is easy to see that the language $L'$ accepted by the modified automaton is the closure of $L$ under free reduction. As $L_0$ is the intersection of $L'$ with the rational set of freely reduced words, $L_0$ is rational. □

## 5. Rational Relations

In addition to the word problem of $G$ introdyced in Section 4 one can consider binary relations induced by the group operation. Let $\Sigma^* \to G$ be a choice of generators, and for each $g \in G$ define a binary relation $\rho_g$ on $\Sigma^*$ by $w\rho_g v$ if $\overline{w}g = \overline{v}$. What can we say about $G$ if $\rho_g$ is a rational subset of $\Sigma^* \times \Sigma^*$? For any $w \in \Sigma^*$, $w\rho_g$ is rational by Corollary 5.5. Since $w\rho_g$ is the word problem of $G$ if $\overline{w} = g^{-1}$, Theorem 4.1 implies that $G$ is finite.

Suppose that instead of demanding that each $\rho_g$ be rational, we require only that the restrictions $\tau_g$ of $\rho_g$ to some fixed rational language $L$ projecting onto $G$ be rational. This line of thought leads to the definition of automatic groups. By Lemma 5.2 $\tau_{g^{-1}} = \tau_g^{-1}$ is rational if $\tau_g$ is, and it follows from Theorem 5.3 that if $\tau_g$ and $\tau_h$ are rational, $\tau_{gh}$ is too. Thus it is sufficient to verify rationality of $\tau_g$ for $g$ in a set of generators of $G$. If we require that for each $g$ in a finite set of generators for $G$, $\tau_g$ is not only rational but is accepted by a special kind of finite automaton over $\Sigma^* \times \Sigma^*$, then $G$ is automatic. See [**7**] for details.

DEFINITION 5.1. Let $M$ and $M'$ be monoids. A relation $\rho : M \to M'$ is rational if its graph is a rational subset of $M \times M'$. A relation between finitely generated free monoids is called a transduction.

From now on we will not distinguish between a relation and its graph. For any relation $\rho : M \to M'$, $m\rho = \{m' \mid m\rho m'\}$; and for $S \subset M$, $S\rho$ is defined similarly. If $m\rho = \{m'\}$, a singleton, we write also $m\rho = m'$.

LEMMA 5.2. *Let $M$ and $M'$ be monoids.*
  (i) *If $\rho : M \to M'$, then the relation $\rho^{-1} : M' \to M$ defined by $m\rho m'$ if and only if $m'\rho' m$ is rational if $\rho$ is.*
  (ii) *If $M = M'$ and $\rho_L = \{(m, m) \mid m \in L\}$ for a rational subset $L \subset M$, then $\rho_L$ is rational.*
  (iii) *If $L$ and $L'$ are rational subsets of $M$ and $M'$ respectively, then $\rho = L \times L'$ is rational.*
  (iv) *If $\rho : M \to M'$ is a partial homomorphism whose domain is a finitely generated submonoid of $M$, then $\rho$ is rational.*

PROOF. To prove the first assertion note that $\rho$ and $\rho'$ correspond under the obvious isomorphism $M \times M' \to M' \times M$, and apply Theorem 2.2. Part (ii) is proved similarly. For part (iii) use the canonical injections $M \to M \times M'$ and $M' \to M \times M'$ together with Theorem 2.2 to conclude that $L \times \{1\}$ and $\{1\} \times L'$ are rational subsets of $M \times M'$. Hence $L \times L' = (L \times \{1\})(\{1\} \times L')$ is rational too. Finally if $\{m_1, \ldots, m_k\}$ generates the domain of $\rho : M \to M'$, then $\rho = \{(m_1, m_1\rho), \ldots, (m_k, m_k\rho)\}^*$. □

THEOREM 5.3. *If $\rho : M \to \Sigma^*$ and $\tau : \Sigma^* \to M'$ are rational relations, so is $\rho \circ \tau$.*

PROOF. Pick sets of generators $A$ and $A'$ for $M$ and $M'$ with both sets of generators containing the identity. Let $\Sigma_1 = \Sigma \cup \{\epsilon\}$. By Theorem 2.6 $\rho$ is accepted by an automaton $\Gamma_\rho$ over $A \times \Sigma_1$ and $\tau$ is accepted by $\Gamma_\tau$ over $\Sigma_1 \times A'$. For each vertex $p$ of $\Gamma_\rho$ add an edge from $p$ to $p$ with label $(1, \epsilon)$ if such an edge is not already present; $\Gamma_\rho$ still accepts $\rho$. Make the corresponding changes to $\Gamma_\tau$. We will see that the automaton $\Gamma$ over $A \times A'$ defined as follows accepts $\rho \circ \tau$.
  (i) The vertices of $\Gamma$ are ordered pairs $(p, q)$ of vertices from $\Gamma_\rho$ and $\Gamma_\tau$ respectively;
  (ii) There is an edge from $(p, q)$ to $(p', q')$ with label $(a, c)$ if and only if for some $b \in \Sigma_1$ there are edges from $p$ to $q$ with label $(a, b)$ and from $p'$ to $q'$ with label $(b, c)$ in $\Gamma_\rho$ and $\Gamma_\tau$ respectively.
  (iii) The initial vertex is $(p_0, q_0)$ where $p_0$ and $q_0$ are the initial vertices of $\Gamma_\rho$ and $\Gamma_\tau$ respectively;
  (iv) The set of terminal vertices is the product of the corresponding sets for $\Gamma_\rho$ and $\Gamma_\tau$.

This construction is illustrated in Figure 3. There we have labelled some single edges with two labels rather than draw two edges with one label each. It is straightforward to show by induction on path length that there is a path in $\Gamma$ from $(p,q)$ to $(p',q')$ with label $(u,w)$ if and only if the following conditions hold.
  (i) For some word $v \in \Sigma^*$ there is a path in $\Gamma_\rho$ from $p$ to $p'$ with label $(u,v)$ and a path in $\Gamma_\tau$ from $q$ to $q'$ with label $(v,w)$.
  (ii) In both path labels $v$ is expressed as a product of elements of $\Sigma_1$ in exactly the same way.

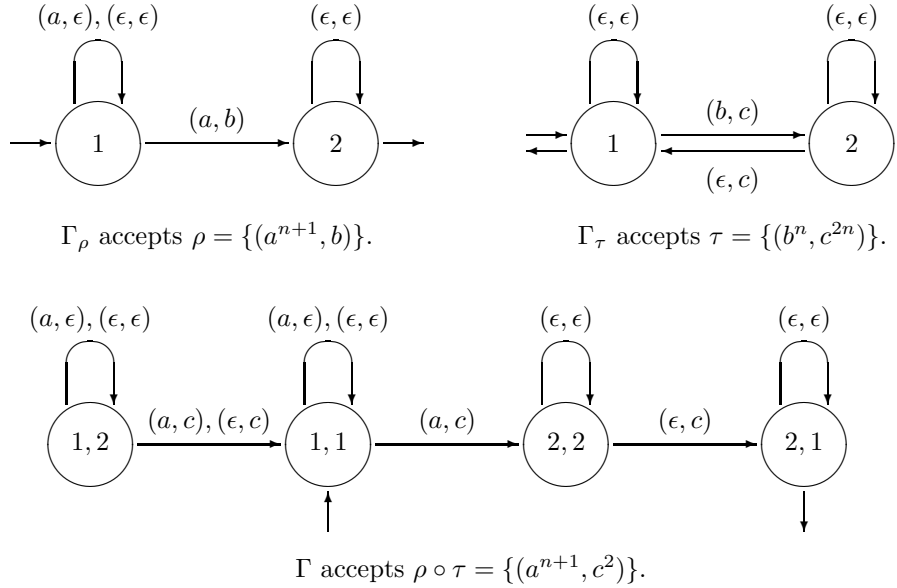

FIGURE 3. Construction of an automaton accepting the composite of two rational relations.

It follows immediately that if $\Gamma$ accepts $(u,v)$, then $u(\rho \circ \tau)v$. For the converse suppose $u\rho v$ and $v\tau w$. There are successful paths in $\Gamma_\rho$ and $\Gamma_\tau$ whose labels are $(u,v)$ and $(v,w)$ respectively. In the label of each path the word $v$ occurs as a product of a sequence of elements of $\Sigma_1$ but the two sequences can differ by insertion and deletion of $\epsilon$'s. Using the extra edges we have added, the successful paths can be adjusted to make the two sequences identical. □

COROLLARY 5.4. *Rational transductions are closed under composition.*

COROLLARY 5.5. *If $\rho : \Sigma^* \to M$ is a rational relation and $L \subset \Sigma^*$ is rational, then $L\rho$ is rational too.*

PROOF. $L\rho$ is the projection of $\rho_L \circ \rho$ from $\Sigma^* \times M$ to $M$. Apply Lemma 5.2 and Theorem 2.2. □

## 6. Families of Languages

In this section we use rational relations to define collections of languages determined by a monoid $M$ and a distinguished subset $X$. For lack of a better word we call these collections families. In the next two sections we see how particular choices of $M$ and $X$ yield context-free languages and stack-automaton languages.

DEFINITION 6.1. The family of languages $\mathcal{F}(M, X)$ determined by a monoid $M$ and a nonempty subset $X \subset M$ is the collection of images $X\rho$ as $\rho$ ranges over all rational relations from $M$ to finitely generated free monoids.

If $M = \{1\}$, then rational relations from $M$ to $\Sigma^*$ are sets of the form $\{1\} \times L$ where $L$ is any rational language over $\Sigma$. It follows that $\mathcal{F}(\{1\}, \{1\})$ is the collection of all rational languages.

By Theorem 2.6 $L \subset \Sigma^*$ is in $\mathcal{F}(M, X)$ if and only if there is an automaton $\Gamma$ over $M \times \Sigma^*$ such that $L$ is the collection of all $w \in \Sigma^*$ with $(m, w)$ the label of a successful path in $\Gamma$ for some $m \in X$. We will say in this case that $\Gamma$ accepts $w$. We have two notions of acceptance. $\Gamma$ accepts both a rational relation $\rho \subset M \times \Sigma^*$ and the language $L = X\rho \subset \Sigma^*$.

We see that the languages in $\mathcal{F}(M, X)$ are characterized concretely in terms of automata at least if $M$ and $X$ are themselves sufficiently tangible. A more abstract characterization is given in the next theorem. When we speak of the image of a language $L \subset \Sigma^*$ under a transduction, we refer to a transduction of the ambient free monoid $\Sigma^*$. A class of languages is closed under rational transduction if whenever $L$ is a language in that class, then the image of $L$ under rational transduction is in the class too. Likewise a class is closed under union if whenever languages $L$ and $L'$ defined over the same alphabet are in the class, so is $L \cup L'$.

THEOREM 6.2. *Let $\mathcal{F} = \mathcal{F}(M, X)$ be a family of languages. $\mathcal{F}$ contains all rational languages and is closed under rational transduction and union. Conversely every collection of languages which contains a nonempty language and is closed under rational transduction and union is a family of languages.*

PROOF. Pick $x \in X$. Every rational language $L$ over $\Sigma$ is the image $X\rho$ under the transduction $\rho = \{x\} \times L$. By Lemma 5.2 $\rho$ is rational, and so $L \in \mathcal{F}$.

To show that $\mathcal{F}$ is closed under union we must show that if $L$ and $L'$ are languages over $\Sigma$ in $\mathcal{F}$, then so is $L \cup L'$. We have $L = X\rho$ and $L' = X\rho'$ for some rational relations $\rho : M \to \Sigma^*$ and $\rho' : M \to \Sigma^*$. It follows that $L \cup L' = X(\rho \cup \rho')$; and since the union of rational relations is a rational relation, we are done. Likewise if $\tau : \Sigma^* \to \Delta^*$ is a rational transduction, then $L\tau = X(\rho \circ \tau) \in \mathcal{F}$ by Theorem 5.3.

For the last part of the theorem let $\mathcal{C}$ be a collection of languages containing a nonempty language and closed under transduction and union. Observe that if we alter a language by substituting new letters for its alphabet, then each version of the language is the image of the other under a homomorphism. Since these homomorphisms are rational transductions by Lemma 5.2, we may pick a subcollection $\mathcal{C}_0$ of languages defined over disjoint alphabets such that every language in $\mathcal{C}$ is the image under homomorphism of some language in $\mathcal{C}_0$.

Now let $M$ be the free monoid over the union of the alphabets of languages in $\mathcal{C}_0$ and define $X$ to be the corresponding union of languages. Since $\mathcal{C}$ contains a nonempty language, $X \neq \phi$. Take $\mathcal{F} = \mathcal{F}(M, X)$. Every language in $\mathcal{C}_0$ is the image of $X$ under the partial identity homomorphism whose domain is the free submonoid of $M$ generated by the alphabet corresponding to that language. Hence every language in $\mathcal{C}$ is the image of $X$ under a partial homomorphism whose domain is a finitely generated free submonoid of $M$. Consequently Lemma 5.2 implies $\mathcal{C} \subset \mathcal{F}$.

It remains to show that $\mathcal{F} \subset \mathcal{C}$. Suppose $L = X\rho$ for some rational relation $\rho : M \to \Sigma^*$. By Theorem 2.2 $\rho \subset M_0 \times \Sigma^*$ where $M_0$ is the free monoid over the union of a finite number of the alphabets used to define $M$. Consequently $L = (X \cap M_0)\rho$. By definition of $X$, $X \cap M_0$ is a finite union of languages in $\mathcal{C}_0$. As $\mathcal{C}$ is closed under union, $X \cap M_0 \in \mathcal{C}$; and as $\mathcal{C}$ is closed under transduction, $(X \cap M_0)\rho \in \mathcal{C}$. □

The closure properties of the preceding theorem resemble some well-known combinations from computer science. If we drop closure under union, we have a full trio; if we add closure under all rational operations, then we obtain the definition of a full abstract family of languages.

COROLLARY 6.3. *Every family of languages is closed under intersection with rational sets.*

PROOF. If $L$ and $R$ are languages over $\Sigma$, then $L \cap R = L\rho_R$ where $\rho_R : \Sigma^* \to \Sigma^*$ is defined by $\rho_R = \{(w, w) \mid w \in R\}$. If $R$ is rational, then by Lemma 5.2 $\rho_R$ is a rational transduction. □

The next theorem shows that it makes sense to speak of the word problem of a group as being in a particular family.

THEOREM 6.4. *Let $\mathcal{F}$ be a family of languages. If the word problem for $G$ with respect to one choice of generators is in $\mathcal{F}$, then the word problem with respect to any choice of generators is also. Further the class of groups whose word problems lie in $\mathcal{F}$ is closed under isomorphism, finitely generated subgroup and finite extension.*

PROOF. Suppose $\sigma : \Sigma^* \to G$ is a choice of generators for $G$ and $\mu : \Delta^* \to H$ is one for the subgroup $H \subset G$. Let $L_G = 1\sigma^{-1}$ and $L_H = 1\mu^{-1}$. There is a monoid homomorphism $f : \Delta^* \to \Sigma^*$ such that $\mu = f \circ \sigma$, and consequently $L_H = L_G f^{-1}$. It follows that $L_G \in \mathcal{F}$ implies $L_H \in \mathcal{F}$. Taking $H = G$, we see that whether or not $L_G \in \mathcal{F}$ depends only on $G$ and not on the choice of generators.

Clearly the class of groups whose word problems lie in $\mathcal{F}$ is closed under isomorphism. Thus it remains only to prove the last assertion. At this point we are free to make any convenient choices of generators for $G$ and $H$. Choose $\sigma : \Sigma^* \to G$ to be a choice of generators with formal inverses. Let $\Gamma$ be the Schreier diagram of $H$ in $G$, i.e., the directed graph whose vertices are the right cosets of $H$ in $G$ and with edges $Hg \xrightarrow{a} Hg'$ whenever $a \in \Sigma$ and $Hg\overline{a} = Hg'$. Taking $H$ to be the initial and single terminal vertex of $\Gamma$, we obtain a deterministic automaton which accepts the language of all words in $\Sigma^*$ representing elements of $H$.

$\Gamma$ admits a spanning tree $\Gamma_0$ rooted at $H$ and with all edges directed away from the root. As in the proof of Theorem 4.2 the edges of $\Gamma$ not in $\Gamma_0$ correspond to generators of $H$. More precisely for each edge $e$ not in $\Gamma_0$ let $w$ be the label of the path in $\Gamma_0$ from $H$ to the source of $e$ and $v$ the path to the target of $e$. If $e$ has label $a$, then $\overline{wa}(\overline{v})^{-1} \in H$, and the collection of these elements as $e$ ranges over all edges not in $\Gamma_0$ generates $H$.

By hypothesis $L_H \in \mathcal{F}$, and we may assume that we have a choice of generators $\mu : \Delta^* \to H$ such that $\Delta$ maps bijectively to the set of generators obtained in the preceding paragraph. That is, the elements of $\Delta$ correspond to the edges in $\Gamma - \Gamma_0$. Alter the edge labels of $\Gamma$ as follows. Replace each label $a$ of an edge in $\Gamma_0$ by $(a, \epsilon)$, and replace each label $a$ of an edge in $\Gamma - \Gamma_0$ by $(a, d)$ where $d$ is the element of $D$ corresponding to that edge. $\Gamma_0$ now accepts a rational transduction $\rho : \Sigma^* \to \Delta^*$ which rewrites each word in $\Sigma^*$ representing an element of $H$ as a word in $\Delta^*$ representing the same element. In fact since $\Gamma$ was originally deterministic, $\rho$ is a partial function. It follows immediately that $w \in L_G$ if and only if $w\rho$ is nonempty and $w\rho \in L_H$. Thus $L_G = L_H \rho^{-1} \in \mathcal{F}$. $\square$

We have seen that groups with rational word problem are finite, and Muller and Schupp [**21**] have shown that groups with context-free word problem are finite extensions of finitely generated free groups. On the other hand many interesting classes of groups are not closed under taking finitely generated subgroups. If we want to study word problems for these classes, we will need to weaken our notion of family.

## 7. Context-Free Languages

In this section we introduce context-free languages. Fix an infinite countable set $\Omega$, and let $\Omega^*$ be the free monoid over $\Omega$. For each word $w \in \Omega^*$ define partial functions $P_w$ and $Q_w$ from $\Omega^*$ to itself by

$$uP_w = uw \text{ for all } u \in \Omega^* \quad uQ_w = v \text{ if } u = vw.$$

We math think of $P_w$ as pushing $w$ onto a stack and $Q_w$ as popping $w$ off a stack. Let $M_{cf}$ be the monoid of partial functions generated by the $P$'s and $Q$'s. $M_{cf}$ has unit element $1 = P_\epsilon = Q_\epsilon$ and zero element 0 (the empty function). The properties listed below are easily verified.

LEMMA 7.1. $P_w P_v = P_{wv}$ and $Q_w Q_v = Q_{vw}$. Further $P_w Q_v = 0$ unless one of the words $w, v$ is a suffix on the other. Every element of $M_{cf}$ except $0$ can be written uniquely as $Q_w P_v$. A product of two of these normal forms equals $0$ except in the following cases.

$$\begin{aligned} Q_w P_{zx} Q_x P_y &= Q_w P_{zy} \\ Q_w P_v Q_{zv} P_y &= Q_{zw} P_y \end{aligned}$$

DEFINITION 7.2. $\mathcal{F}_{cf} = \mathcal{F}(M_{cf}, \{1\})$ is the family of context-free languages.

DEFINITION 7.3. A pushdown automaton is a finite automaton over $M_{cf} \times \Sigma^*$ where $\Sigma^*$ may be any finitely generated free monoid. If $\rho : M_{cf} \to \Sigma^*$ is the rational relation accepted by a pushdown automaton, we say that the automaton accepts the context-free language $1\rho$.

Let us compare Definition 7.3 with the usual definition of a pushdown automaton.

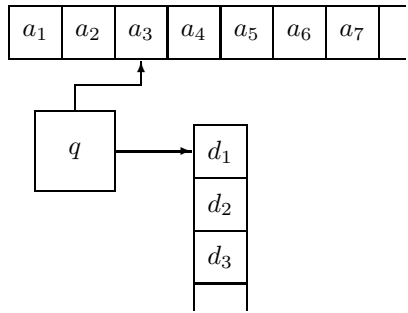

FIGURE 4. A standard pushdown automaton $\mathcal{P}$.

Figure 4 shows a standard pushdown automaton $\mathcal{P}$. $\mathcal{P}$ has a read-only input tape, a pushdown stack, and a finite number of internal states $q$. In Figure 4 the input tape is horizontal and the stack is vertical. At any instant $\mathcal{P}$ is in state $q$ and is scanning a square on its input tape. Based on $q$, the contents of the square being scanned, and the top element of the stack, $\mathcal{P}$ makes a move. During a move $\mathcal{P}$ may pop, i.e., remove, the top element of its stack, push a finite string of symbols onto the stack, move on to the next square of the input tape, and enter a new internal state. It need not do all these things during each move. When a symbol is popped from the stack, the symbol below it becomes the new top of the stack. At different points is its computation $\mathcal{P}$ may have just one possible move, a choice of moves, or no moves at all. In particular $\mathcal{P}$ cannot move if its stack is empty.

At the beginning of a computation the input word is written left justified on the input tape, and the stack is empty except for a special end of stack symbol. $\mathcal{P}$ is placed in a designated initial state and set to scan the first square of its input tape. If there is a sequence of moves which ends with empty stack and $\mathcal{P}$ in one of a number of distinguished final states, and if during this sequence $\mathcal{P}$ reads all its input, then $\mathcal{P}$ accepts the input word. There may be other sequences of moves which end differently or not at all, but that does not matter.

Figure 2 shows a pushdown automaton in our sense; call it $\mathcal{Q}$. Let us view $\mathcal{Q}$ as a standard pushdown automaton and see how it accepts $\{a^n b^n\}$. $\mathcal{Q}$ starts by pushing one $d$ onto the stack for each $a$ it sees in the input. For each succeeding $b$ in the input, $\mathcal{Q}$ tries to pop a $d$ from the stack. If the number of $b$'s equals the number of $a$'s, $\mathcal{Q}$ accepts by emptying the stack and entering its final state. If there are too many $b$'s, $\mathcal{Q}$ empties the stack before reading all its input; and if there are too few, $\mathcal{Q}$ does not empty its stack.

Here is a more formal and shorter demonstration.

THEOREM 7.4. *The pushdown automaton of Figure* 2 *accepts* $L = \{a^n b^n \mid n \geq 0\}$.

PROOF. It is immediate from Figure 2 that the labels of paths from initial vertex to terminal vertex are $\{(P_d^i Q_d^j, a^i b^j)\}$. By Lemma 7.1 $P_d^i Q_d^j = 1$ if and only if $i = j$. □

In the remainder of this section we present some standard facts about context-free languages; but before starting we want to mention a variation on our original definition of pushdown automaton. According to Definition 7.3 a pushdown automaton is an automaton over $M_{cf} \times \Sigma^*$, but by Theorem 2.6 we may restrict ourselves to automata over

$$\bigl(\{P_d, Q_d \mid d \in \Omega\} \cup \{1\}\bigr) \times \bigl(\Sigma \cup \{\epsilon\}\bigr).$$

These automata proceed at a sedate pace, reading at most one input letter and performing at most one stack operation per move. Compared to them our original pushdown automata might be called accelerated.

Now here are the results promised above.

THEOREM 7.5. $\mathcal{F}_{cf}$ *is a full abstract family of languages.*

PROOF. By Theorem 6.2 it suffices to show that if $L$ and $L'$ are context-free languages over $\Sigma$, then so are $LL'$ and $L^*$. Suppose $\Gamma$ and $\Gamma'$ are pushdown automata accepting $L$ and $L'$ respectively. Let $d_1, \ldots, d_n$ be the elements of $\Omega$ which occur in labels of $\Gamma$ and $\Gamma'$, and pick $e \in \Omega$ distinct from the $d_i$'s. Alter $\Gamma$ as indicated in Figure 5 to create a pushdown automaton $\Gamma_1$. More precisely add a new initial vertex joined to the old one by an edge labelled $(P_e, \epsilon)$, and join all old terminal vertices to a single new terminal vertex by edges labelled $(Q_e, \epsilon)$. Clearly $\Gamma_1$ accepts $(P_e Q_v P_w Q_e, u)$ if and only if $\Gamma$ accepts $(Q_v P_w, u)$. It follows from Lemma 7.1 and our choice of $e$ that $\Gamma_1$ accepts $(1, u)$ for all $u \in L$ and that everything else accepted has the form $(0, y)$ for some $y \in \Sigma^*$. Modify $\Gamma'$ similarly to obtain $\Gamma'_1$ which accepts $(1, u)$ for all $u \in L'$ and possibly $(0, y)$ with $y \in \Sigma^* - L'$. These pushdown automata may be used as indicated in Figure 6 to construct pushdown automata accepting $LL'$ and $L^*$. □

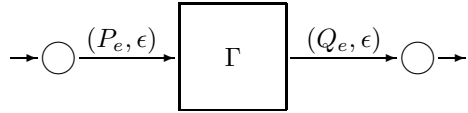

FIGURE 5. Construction of the pushdown automaton $\Gamma_1$.

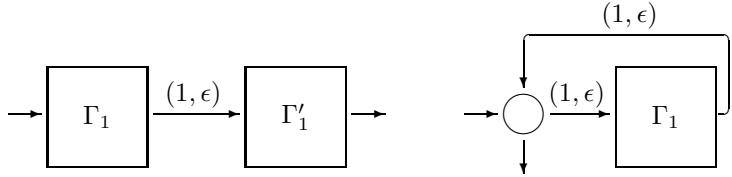

FIGURE 6. Construction of pushdown automata accepting $LL'$ and $L^*$.

Our first impulse in proving Theorem 7.5 is to connect the unmodified automata $\Gamma$ and $\Gamma'$ the way $\Gamma_1$ and $\Gamma'_1$ are connected in Figure 6. However, it might be that $\Gamma$ accepts $(P_w, u)$ and $\Gamma'$ accepts $(Q_w, v)$. In this case $\Gamma$ followed by $\Gamma'$ would accept $uv$ not necessarily in $LL'$. Thus modifications are necessary.

The following theorem is a pumping lemma for context-free languages. Its usual first application is to show that the language $\{a^n b^n c^n\}$ is not context-free.

THEOREM 7.6. *For every context-free language $L$ there is an integer $k$ such that any word in $L$ of length greater than $k$ can be written as $uvwxy$ so that*
  (i) $v \neq \epsilon$ or $x \neq \epsilon$;
  (ii) $|vwx| \leq k$;
  (iii) $uv^i wx^i y \in L$, for all $i \geq 0$.

It is an immediate consequence of the Pumping Lemma that a group whose word problem is context-free admits a Dehn's Algorithm solution to its word problem. Thus such groups are word hyperbolic. The work of Muller and Schupp shows that in fact they are virtually free.

Let $L$ be a context-free language accepted by a pushdown automaton $\Gamma$ over $M_{cf} \times \Sigma^*$. We will prove the Pumping Lemma by showing that a sufficiently long successful path in $\Gamma$ must have one of two structures. Either for some $s \in \Omega^*$ the path contains a loop with label $(P_s, v)$ followed by a segment with label $(1, w)$ followed by a loop with label $(Q_s, x)$; or the path contains a loop with label $(1, v)$. The latter case corresponds to the possibility $w = x = \epsilon$.

Since $M_{cf}$ is generated by $\{P_d, Q_d \mid d \in \Omega\} \times (\Sigma \cup \{\epsilon\})$, we may assume that $\Gamma$ is defined over this set. Our proof depends on an analysis of the word problem of $M_{cf}$ with respect to the generating set $\{P_d, Q_d \mid d \in \Omega\}$.

LEMMA 7.7. *Let $A_1 \ldots A_n$ be a nonempty word in the generators $\{P_d, Q_d \mid d \in \Omega\}$ representing $1 \in M_{cf}$. One of the following conditions holds.*
  (i) *$A_1 = P_d$, $A_n = Q_d$ for some $d \in \Omega$, and either $n = 2$ or $n > 2$ and $A_2 \cdots A_{n-1}$ represents 1; or*
  (ii) *For some $i$ with $1 < i < n$, $A_1 \cdots A_i$ and $A_{i+1} \cdots A_n$ both represent 1.*
*Further either $n = 2$ or $A_1 \ldots A_n$ has a proper subword of length at least $n/2$ which represents 1.*

PROOF. It follows from Lemma 7.1 that $A_1 \ldots A_n$ contains a subword $P_d Q_d$. With this observation the first part of the following lemma is easily proved by induction on $n$. The first part implies the last part. □

Now let $c$ be the number of vertices of $\Gamma$ and choose a positive integer contant $k > 8c$. We permit ourselve to increase $k$ as we proceed. Consider a word $z \in L$ of length at least $k$, and pick a successful path $\gamma$ in $\Gamma$ of minimal length with label $(1, z)$. It follows from Lemma 7.7 that $\gamma$ has a subpath $\gamma_0$ of length between $k$ and $k/2$ and with label $(1, z_0)$ for some subword $z_0$ of $z$. We will analyze this subpath.

First suppose $\gamma_0$ has a sequence of $c$ consecutive nonempty subpaths $\gamma_1, \ldots, \gamma_c$ with labels $(1, z_i)$ for $1 \leq i \leq c$. Let path $\gamma_i$ extend from vertex $p_{i-1}$ to vertex $p_i$. For some $i, j$ with $0 \leq i < j \leq c$, the vertices $p_i$ and $p_j$ must be equal. Consequently $\gamma$ has a nonempty loop of length at most $k$ with label $(1, v)$ for some $v \in \Sigma^*$. If we alter $\gamma$ by traversing the loop $i \geq 0$ times, we still have a successful path. In particular $v \neq \epsilon$ lest there be a shorter successful path than $\gamma$ with label $(1, z)$. Thus we have $z = uvy$, $v \neq \epsilon$, $|v| \leq k$ and $uv^i y \in L$ for all $i \geq 0$.

Now assume $\gamma_0$ does not have a sequence of $c$ consecutive subpaths as above. Write $\gamma_0$ as a sequence of as many nonempty subpaths as possible such that each one has label $(1, z')$ for some $z' \in \Sigma^*$. Since there are fewer than $c$ of these subpaths, one of them, call it $\gamma'_O$, has length at least $k/(2c)$. Further $\gamma'_0$ cannot itself be decomposed as a product of nonempty subpaths with labels $(1, z')$. By Lemma 7.7 the sequence of lefthand components of edge labels of $\gamma'_0$ must have the form $P_d B_1 \ldots B_m Q_d$ for some nonempty word $B_1 \ldots B_m$ in the generators $\{P_d, Q_d \mid d \in \Omega\}$ representing $1 \in M_{cf}$. The proper subpath $\gamma_1$ of $\gamma'_0$ corresponding to $B_1 \ldots B_m$ has label $(1, z_1)$ and by assumption cannot be decomposed as a sequence of $c$ consecutive subpaths satisfying the conditions of the last paragraph.

We can repeat the argument of the last paragraph with $\gamma_1$ in place of $\gamma_0$ to obtain $\gamma_2$. It is clear that if $k$ is large enough, we can find paths $\gamma_0, \gamma_1, \ldots, \gamma_{c^2}$ each a proper subpath of the preceeding path and each with label $(1, z')$ for some $z' \in \Sigma^*$. Since there are only $c^2$ pairs of vertices from $\Gamma$, two of these paths must begin at the same vertex and end at the same vertex. It follows that for some words $r, s, r's' \in \Omega^*$, $\gamma$ has a subpath consisting of a loop with label $(Q_r P_s, v)$ followed by a segment with label $(1, w)$ followed by a loop with label $(Q_{r'} P_{s'}, x)$.

Further at least one of the loops is nonempty and $Q_r P_s Q_{r'} P_{s'} = 1$. The latter condition implies the first loop has label $P_s$ and the second $Q_s$. Consequently if we alter $\gamma$ by traversing each loop $i$ times instead of one, we we still have a successful path. In particular $vx \neq \epsilon$ lest there be a shorter successful path than $\gamma$ with label $(1, z)$. Thus we have $z = uvwxy$, $vx \neq \epsilon$, $|vwx| \leq k$ and $uv^i wx^i y \in L$ for all $i \geq 0$.

Our argument has shown that if $k$ is sufficiently large, all the conditions of the pumping lemma hold.

## 8. Stack Automata

In this section we briefly consider stack automata. A standard stack automaton is a pushdown automaton which can read but not alter the contents of its stack below the top. Figure 7 shows a stack automaton reading an element from the interior of its stack. The stack is shown in two pieces in order to suggest how two words may be used to represent both the contents of the stack and the stack location being read by the automaton.

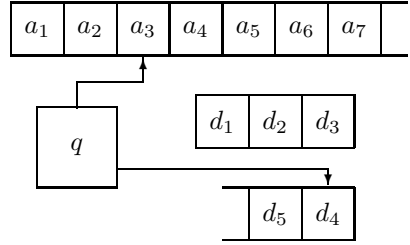

Figure 7. A standard stack automaton.

The appropriate monoid for a formal definition of a stack automaton is built out of the monoid $M_{cf}$ of partial functions $P_w$ and $Q_w$ defined in the preceding section. Since we will need to know when a stack is empty, we add a new partial function

$$uE = \epsilon \text{ if } u = \epsilon$$

and define $M_1$ to be the submonoid of partial functions generated by $M_{cf} \cup \{E\}$.

Lemma 8.1. *Every element of $M_1$ except $0$ can be written as $Q_w P_v$ or $Q_w E P_v$.*

*A product of two such terms equals $0$ except in the following cases.*

$$\begin{aligned} Q_w P_{zx} Q_x P_y &= Q_w P_{zy} \\ Q_w P_v Q_{zv} P_y &= Q_{zw} P_y \\ Q_w E P_{zx} Q_x P_y &= Q_w E P_{zy} \\ Q_w P_v Q_{zv} E P_y &= Q_{zw} E P_y \\ Q_w E P_x Q_x E P_y &= Q_w E P_y \end{aligned}$$

Consider the monoid $M_1 \times M_1$ acting from the right on $\Omega^* \times \Omega^*$. We think of $(w,v) \in \Omega^*$ as representing a stack which contains $wv^r$, where $v^r$ is $v$ read from right to left. In Figure 7 $w = d_1 d_2 d_3$, and $v = \ldots d_5 d_4$. The leftmost symbol of $w$ is at the top of the stack, and the automaton is reading the rightmost symbol of $v$ inside the stack. The operations of moving up and down in the stack are $(Q_x, P_x)$ and $(P_x, Q_x)$ respectively. Ordinary stack operations are permitted only at the top of the stack, that is when $w = \epsilon$. They are $(E, P_x)$ and $(E, Q_y)$. $M_{sa}$ is the submonoid of $M_1 \times M_1$ generated by all these elements.

DEFINITION 8.2. $\mathcal{F}_{sa} = \mathcal{F}(M_{sa}, \{(E,1)\})$ is the family of stack-automaton languages.

In other words the stack-automaton languages are the images of $(E,1)$ under rational relations $\rho : M_{sa} \to \Sigma^*$.

DEFINITION 8.3. A stack automaton is a finite automaton over $M_{sa} \times \Sigma^*$ where $\Sigma^*$ may be any finitely generated free monoid. If $\rho : M_{sa} \to \Sigma^*$ is the rational relation accepted by a stack automaton, we say that the automaton accepts the stack-automaton language $(E,1)\rho$.

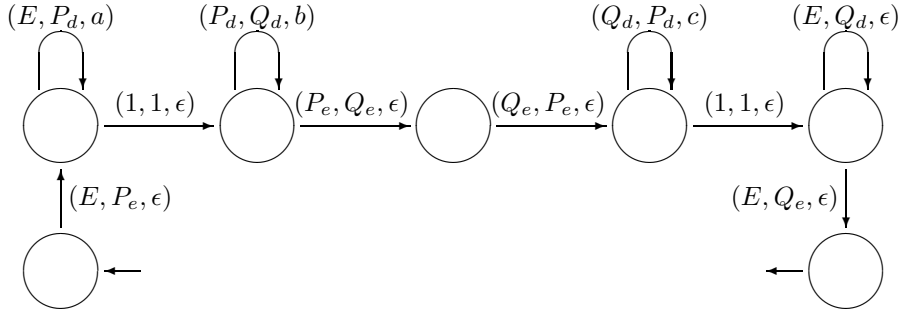

FIGURE 8. A stack automaton $\mathcal{S}$ accepting $\{a^n b^n c^n\}$.

Figure 8 shows a stack automaton $\mathcal{S}$ which accepts the language $\{a^n b^n c^n\}$. $\mathcal{S}$ works by
  (i) Pushing an end of stack marker onto its stack.
  (ii) Pushing one $d$ onto its stack for each $a$;

(iii) Matching $d$'s on its stack against successive $b$'s until it reaches the bottom of the stack;
(iv) Matching $d$'s on its stack against successive $c$'s until it reaches the top of the stack;
(v) Emptying the stack.

LEMMA 8.4. *The automaton $\mathcal{S}$ in Figure* 8 *accepts* $\{a^n b^n c^n\}$.

PROOF. Any successful path has label

$$(E^{i+1}P_d^j P_e Q_e Q_d^k E^{m+1}, P_e P_d^i Q_d^j Q_e P_e P_d^k Q_d^m Q_e, a^i b^j c^k,) = $$
$$(E P_d^j Q_d^k E, P_e P_d^i Q_d^j Q_e P_e P_d^k Q_d^m Q_e, a^i b^j c^k).$$

The middle entry is $E$ if and only if $j = k$, and last entry reduces to 1 exactly when $i = j$ and $k = m$. □

## 9. Grammars

In this section we briefly discuss how grammars are used to generate languages. We will show that the languages we have called context-free are exactly the languages generated by context-free grammars and thus coincide with the standard context-free languages.

DEFINITION 9.1. A grammar consists of
(i) A finite set of terminal symbols $\Sigma$;
(ii) A disjoint finite set of nonterminals $\Delta$;
(iii) A start symbol $S \in \Delta$;
(iv) A finite set of productions $\mathcal{P} = \{(\alpha, \beta)\} \subset (\Sigma \cup \Delta)^* \times (\Sigma \cup \Delta)^*$ such that each $\alpha$ contains at least one symbol from $\Delta$.

Definition 9.1 is the usual one, but many variations exist. $(\Sigma \cup \Delta)^*$ is called the set of *sentential forms*. Each grammar determines a binary relation on sentential forms. This relation is defined in stages. First we say that a sentential form $\gamma$ *directly derives* another sentential form $\gamma'$ with respect to a given grammar if $\gamma'$ is obtained by substituting the righthand side of some production for the lefthand side in $\gamma$. The reflexive transitive closure of direct derivation is called *derivation*. We write $\gamma \to \gamma'$ for direct derivation and $\gamma_1 \xrightarrow{*} \gamma_2$ for derivation. The language generated by a grammar is the set of words over the terminal alphabet derivable from the start symbol. From now on we will write $\alpha \to \beta$ for productions instead of $(\alpha, \beta)$.

DEFINITION 9.2. A grammar is regular if its productions are of the form $A \to aB$ or $A \to a$ where $A$ and $B$ are nonterminals and $a$ is a terminal.

Regular grammars are just another way of looking at finite automata over free monoids, and they generate exactly the regular languages. Indeed suppose that $\Gamma$ is an automaton over $\Sigma$ accepting $L \subset \Sigma^*$. Assume its vertices are $p, q, \ldots$, and label them with distinct symbols $P, Q, \ldots$ not in $\Sigma$. Let $P_0$ be the label of the initial vertex $p_0$. To each edge $p \xrightarrow{a} q$ associate the production $P \to aQ$. If $q$ is a terminal state, add the production $P \to a$. A successful path

$$p_0 \xrightarrow{a_1} p_1 \xrightarrow{a_2} \cdots \xrightarrow{a_{m-1}} p_{m-1} \xrightarrow{a_m} p_m$$

to a terminal vertex $p_m$ is essentially the same as a derivation

$$P_0 \to a_1 P_2 \to a_1 a_2 P_2 \to \cdots \to a_1 \ldots a_{m-1} P_{m-1} \to a_1 \ldots a_m.$$

Thus the grammar consisting of
   (i) Terminal alphabet $\Sigma$;
   (ii) Nonterminal alphabet $\{P, Q, \ldots\}$;
   (iii) Start symbol $P_0$;
   (iv) Productions $\{P \to aQ\} \cup \{P \to a\}$ as above

generates $L$. This argument can be reversed to show that any regular grammar can be converted to a finite automaton accepting the language generated by the grammar. If the grammar has a production $P \to a$ but no corresponding production $P \to aQ$, one can chose a new nonterminal $Q$ and add the missing production. Since there are no productions with $Q$ on the lefthand side, the language generated is not changed.

Now we proceed to context-free grammars.

DEFINITION 9.3. A grammar is context-free if the lefthand side of every production consists of a single nonterminal.

Let $F$ be a free group of rank two and $W(F)$ the word problem of $F$ with respect to a choice of generators $\Sigma = \{a, a^{-1}, b, b^{-1}\}$, $\Sigma^* \to F$ such that $\overline{a}$ and $\overline{b}$ generate $F$ freely. We will show that $W(F)$ is generated by a context-free grammar. It follows from this result and Theorem 6.4 that every finitely generated group with a free subgroup of finite index has a context-free word problem.

The key to proving that $W(F)$ is context-free is its recursive nature. Any nonempty word which projects to the identity is either a product of shorter words which project to the identity, or has one of the forms $ava^{-1}, bvb^{-1}, a^{-1}va, b^{-1}vb$ for some $v \in W(F)$. Accordingly we let $\mathcal{G}$ be the grammar with productions

$$S \to \epsilon, \quad S \to SS, \quad S \to aSa^{-1}, \quad S \to bSb^{-1}, \quad S \to a^{-1}Sa, \quad S \to b^{-1}Sb.$$

A sample derivation $S \xrightarrow{*} abb^{-1}b^{-1}ba^{-1}$ with respect to $\mathcal{G}$ is given by

$$S \to aSa^{-1} \to aSSa^{-1} \to abSb^{-1}Sa^{-1} \to abb^{-1}Sa^{-1}$$
$$\to abb^{-1}b^{-1}Sba^{-1} \to abb^{-1}b^{-1}ba^{-1}.$$

To show that $\mathcal{G}$ is indeed a grammar generating $W(F)$ there are two steps, both easy. First prove by induction on derivation length that if $S \xrightarrow{*} w$, then $\overline{w} = 1$. Second use the observations on $W(F)$ made above to prove by induction on word length that if $w \in W(F)$, then there is a derivation $S \xrightarrow{*} w$.

Our next goal is to show that Definition 7.2 defines precisely the class of languages generated by context-free grammars.

LEMMA 9.4. *If $L$ is generated by a context-free grammar, then $L$ is generated by a context-free grammar in which the righthand side of every production consists entirely of nonterminals or is a single terminal or is $\epsilon$.*

PROOF. For each terminal $a$ add a new variable $A$ and a new production $A \to a$ to the grammar. Replace all occurrences of terminals in the original productions by the corresponding new variables. It is straightforward to check that $L$ remains unchanged. □

LEMMA 9.5. *If there is a derivation $S \overset{*}{\to} w$ with respect to a context-free grammar $\mathcal{G}$, then there is a derivation in which productions are always applied to the leftmost nonterminal.*

PROOF. We show by induction on derivation length that for all nonterminals $A$ and derivations $A \overset{*}{\to} w$ of words in the terminal alphabet there is a leftmost derivation. Clearly any derivation of length one is leftmost for it must consist of the single step $A \to w$. The first step of any longer derivation $A \overset{*}{\to} w$ must be $A \to w_0 A_1 w_1 \cdots A_n w_n$ where the $A_i$'s are nonterminals and the $w_i$'s are words in the terminal alphabet. The remaining steps in the derivation amount to derivations $A_i \overset{*}{\to} u_i$ such that $w = w_0 u_1 \cdots u_n w_n$. By induction on derivation length, these derivations may be taken to be leftmost from which it is clear that there is a leftmost derivation $A \overset{*}{\to} w$. □

THEOREM 9.6. *A formal language is generated by a context-free grammar if and only if it is accepted by a pushdown automaton.*

PROOF. First suppose that $L \subset \Sigma^*$ is generated by a context-free grammar $\mathcal{G}$; we will use $\mathcal{G}$ to construct a pushdown automaton accepting $L$. There is no loss of generality in assuming that $\mathcal{G}$ satisfies the condition of Lemma 9.4 and that its nonterminals are in the set $\Omega$ introduced at the beginning of Section 7. For any sentential form $\alpha$, $\alpha^r$ stands for the reverse of $\alpha$, i.e., $\alpha$ read from right to left. Define $\Gamma$ to be a pushdown automaton with an initial state, $q_0$, and a terminal state $q_t$. Add an edge from $q_t$ to $q_t$ with label
  (i) $(Q_A P_{\alpha^r}, \epsilon)$ for each production $A \to \alpha$ in $\mathcal{G}$ with $\alpha$ consisting of nonterminals;
  (ii) $(Q_A, a)$ for each production $A \to a$ with $a \in \Sigma \cup \{\epsilon\}$.
Add one more edge from $q_0$ to $q_t$ with label $(P_S, \epsilon)$ where $S$ is the start symbol for the grammar. Prove by induction on path length that if there is a path from $q_0$ to $q_t$ with label $(P_{\gamma^r}, w)$ then there is a leftmost derivation $S \overset{*}{\to} w\gamma$; and prove the converse by induction on derivation length. In particular there is a successful path with label $(P_\epsilon, w) = (1, w)$ if and only if there is a leftmost derivation $S \overset{*}{\to} w$. It follows from Lemma 9.5 that $L$ is the language accepted by $\Gamma$.

Now suppose $L \subset \Sigma^*$ is a context-free language. By Theorem 2.6 $L$ is accepted by a pushdown automaton $\Gamma$ over $\{P_d, Q_d \mid d \in \Omega\} \times (\Sigma \cup \{\epsilon\})$. Let $p_0$ be the initial vertex of $\Gamma$, and for each pair of vertices $p, q$ define $L_{p,q}$ to be the language accepted by $\Gamma$ with its start state changed to $p$ and with single terminal state $q$. We will construct context-free grammars $\mathcal{G}_{p,q}$ for these languages; the grammars will be identical except for their start symbols. As $L$ is the union over all terminal vertices $q$ of the languages $L_{p_0,q}$, we obtain a grammar for $L$ by adding a new start symbol $S$ and productions $S \to A_{p_0,q}$ to the union of the grammars $\mathcal{G}_{p,q}$ over all terminal vertices $q$.

To define $\mathcal{G}_{p,q}$ start with terminals $\Sigma$ and one nonterminal $A_{p,q}$ for each pair of vertices $p, q$ of $\Gamma$ such that there is a path in $\Gamma$ from $p$ to $q$ with label $(1, w)$ for some $w \in \Sigma^*$. In particular there are nonterminals $A_{p,p}$ for all vertices $p$ of $\Gamma$. Add

(i) Productions $A_{p,q} \to A_{p,r} A_{r,q}$ for all triples $p, r, q$ such that the corresponding nonterminals exist;
(ii) Productions $A_{p,q} \to a A_{r,s} b$ for all nonterminals $A_{p,q}$ such that there is a nonterminal $A_{r,s}$, an edge labeled $(a, P_d)$ from $p$ to $r$, and one labeled $(b, Q_d)$ from $s$ to $q$;
(iii) Productions $A_{p,p} \to \epsilon$ for all vertices $p$.

A straightforward argument using Lemma 7.7 and induction on path length shows that if $(1, w)$ is the label of a path from $p$ to $q$ in $\Gamma$, then there is a derivation $A_{p,q} \stackrel{*}{\to} w$. It is also straightforward to show by induction on derivation length that $A_{p,q} \stackrel{*}{\to} w$ implies that $(1, w)$ is the label of a path from $p$ to $q$ in $\Gamma$. It follows from these results that the grammar $G_{p,q}$ with start symbol $A_{p,q}$ and other components as above generates $L_{p,q}$. □

## 10. Literature

The first paragraph of the introduction mentions several applications of formal language theory to group theory; here are the corresponding citations. Automatic groups and word hyperbolic groups were introduced by Epstein, Cannon, Holt, Levy, Paterson and Thurston [7], and Gromov [14] respectively. Sims' book [24] contains many applications of rational languages to computational group theory, and in Chapter 4 there is an automata-theoretic analysis of finitely generated subgroups of free products of cyclic groups. Muller and Schupp [21] showed that a group has context-free word problem if and only if its word problem is virtually free. For a survey of further developments along these lines including the characterization of virtually cyclic groups by one-counter languages see the survey by Herbst and Thomas [16]. Bridson and Gilman [5] have characterized fundamental groups of geometrizable three-manifolds in terms of indexed languages.

One of the earliest investigations of word problems of groups in terms of formal languages is by Anisimov [1]. Theorem 4.2 is usually credited to Anisimov and Seifert [2]; a special case which extends immediately to the general result is given in Chapter 6 of Johansen's dissertation [18].

There is another concept of automatic group which predates the one referred to above, namely a group of invertible rational transductions of a free monoid. Recall that by Corollary 5.4 rational transductions are closed under composition. Automatic groups in this sense with intermediate growth have been discovered by Grigorchuk [13]. We note that the class of groups whose elements are rational transductions permuting a fixed rational language includes both types of automatic group.

The technique used in this article of defining different kinds of automata in terms of rational subsets of monoids is well-known and probably grew out of the observation that popping an element off a stack is a right inverse to pushing it on. Similar tacks are taken by Brainerd and Landweber [3, Chapter 4], Eilenberg [9, Volume A, Chapter X], Floyd and Biegel [10], Goldstine [12], and Salomaa, Wood and Yu [25]. Most of these cited works include diagrams resembling Figure 2. For standard presentations of language theory one can consult Harrison [15] or Hopcroft and Ullman [17]. Revesz [23] has a slightly different perspective, and Conway [6] gives an extensive treatment of regular languages.

Although we have not mentioned them here, there are many applications of language theory to semigroups and monoids. See Lallement [19], Book and Otto [4], Madlener and Otto [20], as well as Eilenberg's two volumes.

Department of Mathematics, Stevens Institute of Technology, Hoboken, New Jersey

*E-mail address*: rgilman@vaxc.stevens-tech.edu